\documentclass{amsart}
\usepackage{latexsym}                  
\usepackage{exscale}                   

\theoremstyle{plain}

\newtheorem*{lem}{Lemma}
\newtheorem*{cor}{Corollary}
\newtheorem*{lclt}{Local Central Limit Theorem}

\newcommand{\ane}[1]{`#1'}          

\newcommand{\ssst}{\scriptscriptstyle}

\newcommand{\field}[1]{\mathbb{#1}}
\newcommand{\C}{\field{C}}
\newcommand{\N}{\field{N}}
\newcommand{\No}{\N_0}

\newcommand{\Rd}{\field{R}^d}

\newcommand{\Zd}{\field{Z}^d}

\newcommand{\deff}{\stackrel{\rm def}{=}}

\DeclareMathOperator{\Id}{Id}
\DeclareMathOperator{\inter}{int}

\newcommand{\cons}{C}
\newcommand{\const}{C}
\newcommand{\con}{L}
\newcommand{\conste}{K}
\newcommand{\constz}{{K'}}
\newcommand{\klconst}{\delta}
\newcommand{\klconste}{\delta'}
\newcommand{\step}{\ell}

\newcommand{\E}{E}                 
\newcommand{\V}{\Gamma}            
\newcommand{\klew}{\gamma}         
\newcommand{\Nullum}{{\mathcal E}} 

\newcommand{\abs}[1]{\lvert#1\rvert}

\newcommand{\babs}[1]{\bigl\lvert#1\bigr\rvert}

\newcommand{\Babs}[1]{\Bigl\lvert#1\Bigr\rvert}

\newcommand{\babbs}[1]{\biggl\lvert#1\biggr\rvert}

\newcommand{\AG}{G}
\newcommand{\AH}{H}

\newcommand{\huth}{\widehat{\AH}}

\newcommand{\sumx}[1][x]{\sum_{#1\in\Zd}}

\begin{document}
\title{Good Local Bounds for Simple Random Walks}

\author{Christine Ritzmann}
\date{\today\\[1mm]  Angewandte
Mathematik, Universit\"at Z\"urich, Winterthurer
Strasse 190, CH-8057 Z\"urich, Switzerland.
Tel.\,xx41/1/635~58~54, \texttt{chritz@amath.unizh.ch}}

\begin{abstract}
We give a local central limit theorem
for simple random walks on $\Zd$, including Gaussian error estimates.
The detailed proof combines standard large deviation
techniques with  Cram\'er-Edgeworth expansions for lattice distributions.
\end{abstract}
\maketitle

A simple random walk on $\Zd$ is a sequence of independent
identically distributed random variables (steps) on $\Zd$. We here suppose
that the steplength is bounded, that is, the measure describing a single
step has bounded support, and our goal is to give local approximations for
the sum of the first
$n$ steps of the walk. Of course the Central Limit Theorem applies, and
there also exist local CLTs, see
for example Bhattacharya/Rao, \cite{BR}.

In our work on local approximations for polymer measures
(see the joint paper with E.\,Bolt\-hausen \cite{Bori}), we considered
weakly self-avoiding walks as perturbations of simple walks. In order to
obtain good bounds for the perturbations we thus needed good bounds for
the simple walks as well. In particular, we needed Gaussian tail
estimates, which were not provided by \cite{BR}. Therefore we proved a
version of a local CLT, which was tailored for our application there. 

It turns out that such a LCLT is useful not just in the treatment of
polymer measures, but in all cases where
similar perturbation methods are used and hence good local bounds for
simple random walks are needed. 

Therefore we present a local central limit theorem with Gaussian tail
estimates for walks with bounded steplength. This restriction may exclude
some possible applications, but the setting is sufficient for many
purposes and it allows a precise control on the error terms.  A major
advantage over \cite{BR} is that we make the dependence of the constants
explicit.

In order to formulate the theorem, we need some
definitions and notations:
Let $\AG $ be the single step distribution of a random walk on $\Zd$. The
walk has \emph{bounded steplength $\step$} if $\AG(x)=0$ for all
$\abs{x}>\step$. 
The walk is called \emph{maximal}, if the
support of $\AG$ is not contained in any affine hyperplane of $\Zd$.
The walk is called \emph{aperiodic}, if the greatest common divisor of all
$n\in\N$ with $\AG^{*n}(0)>0$ is one.

In the whole paper, $\varphi_\V$ denotes the density of the
centered normal distribution on
$\Rd$ with covariance matrix $\V$, that is,
\begin{equation*}
\varphi_\V (x) = (2\pi)^{-d/2}(\det \V)^{-1/2}
\,\exp[- x\cdot(\V^{-1}\cdot x)/2].
\end{equation*}

Furthermore we abbreviate $\varphi_{(\eta\cdot \Id_d)}$ by
$\varphi_\eta$, that is
\begin{equation*}
\varphi_\eta (x) = (2\pi\eta)^{-d/2} \,\exp\bigl[-x^2/(2\eta)].
\end{equation*}

Note that for $x\in\Rd$ we abbreviate $\abs{x}^k$ by $x^k$.

\begin{lclt}                                      
Let $\AG $ be the single step distribution of a maximal and
aperiodic random walk on $\Zd$ with mean zero and bounded steplength
$\step\in\N$. Let $\V$ denote the covariance matrix of $\AG $.
Then there exist polynomials $P_3$ and $P_6$ of degree three and six,
respectively, and a positive constant
$\const$ such that  for all $x$ in $\Zd$, $\alpha\in (0,1/2)$
and for all  natural $n$,
\begin{align*}                                      
\babs{\AG ^{*n}(x)-\bigl[1+n^{-1/2}P_3 (x/\sqrt{n})
+ n^{-1}P_6 (x/\sqrt{n})\bigr]\varphi_{n\V}(x)} &\le \const
n^{-1-\alpha}\varphi_{n 2d\step^2}(x).
\end{align*}
The coefficients of the polynomials depend (polynomially) on $\V^{-1}$
and on the moments of
$\AG$ up to order four, whereas $\const$ can be chosen
independently of the specific law of $\AG $, depending only on $d$,
$\step$, $\alpha$ and on a lower bound for the smallest eigenvalue of
$\,\V$.
\end{lclt}

Before proving the theorem, we make some comments on the
assumptions. The fact that the measure has bounded support is essential
for the proof and can not be given up easily. The other assumptions are
only technical: By shifting the measure we can arrange zero expectation
and by restricting to a subspace we achieve maximality. For periodic
measures we have to adapt the inversion formula for lattice measures
\eqref{IF1}. The effect of the period is that for fixed $n$ we have non
zero measure of $\AG^{*n}$ only on a sublattice of $\Zd$, which leads to
an additional factor in the approximating term but leaves the rest of
the proof unchanged.

The proof of the LCLT is split in two parts. The first part contains
the large deviation arguments. In the second part we prove the
Cram\'er-Edgeworth expansions, which are used in the first part to
approximate the measures obtained by tilting the original distribution.

\begin{proof}[Proof of the LCLT]
Let
\begin{align*}
Z(t) &\deff \sumx \exp[t\cdot x]\AG (x)\\
\text{and}\qquad I(\xi)&\deff \sup_{t\in\Rd}\{t\cdot \xi - \log Z(t)\}.
\end{align*}

Standard large
deviation theory (see for example \cite{E}) yields a large deviation principle
with entropy function $I$ for the laws of $\AG ^{*n}(x/n)$.
Let $S_\AG$ denote the convex closure of the set of points with nonzero
$\AG$ measure. Then  $I$ is convex on $\Rd$ and even strictly convex on
$\inter S_\AG $, that is, the interior of $S_\AG$. Outside
$S_\AG $, $I$ equals $+ \infty$.

The function $t\mapsto
D\log Z(t)$ is an analytic diffeomorphism from $\Rd$ onto $\inter
S_\AG $ (for a proof see \cite{E}, page 261).
Therefore, for any
$\xi\in \inter S_\AG $, there exists a unique $t_\xi\in\Rd$ with $D\log
Z(t_\xi)=\xi$. Clearly $D \log Z(0)=0$ and $D^2 \log
Z(0)=\V$. For
$\xi \in \inter S_\AG $, we have $I(\xi)= t_\xi\cdot
\xi - \log Z(t_\xi)$.
Now for $t \in \Rd$ denote by $\AG _t$ the tilted measure
\begin{equation*}
\AG _t (x) \deff \frac{\AG (x) \exp[t\cdot x]}{Z(t)}.
\end{equation*}

Using this, we see that for $\xi\deff x/n \in \inter S_\AG $, we can write
\begin{equation}                                            \label{AG-n}
\AG ^{*n} (x) = \exp[-n I(\xi)] \cdot \AG _{t_\xi}^{*n}(x).
\end{equation}

\vspace{1mm}
{\em Case $\abs{x}\le n^{(5-\alpha)/9}$:}

Since $\AG $ is centered and maximal, the
boundary of $S_\AG $ is bounded away from zero. The covariance matrix
$\V_\xi$ of $\AG _{t_\xi}$ depends
analytically on $\xi$ with $\V_0=\V$. Let $\klew$ be the smallest
eigenvalue of $\V$.
We can find a constant $\klconst>0$ such that for all $\xi\in \Rd$ with
$\abs{\xi}\le\klconst$ we have
\begin{itemize}
\item
$\xi \in \inter S_\AG$ and
\item
the smallest eigenvalue of $\V_\xi$ is greater or equal $\klew/2$.
\end{itemize}

We will now argue that $\klconst$ can be chosen depending only on $\step$
and $\klew$, but not on the specific measure $\AG$. There is only a finite
number of possible supports for a maximal and centered measure $\AH$ with
steplength $\step$ on $\Zd$, and the convex hulls of these supports
all contain zero in their interior. Therefore we find a constant
$\klconste$ only depending on $\step$, such that for all
$\abs{\xi}\le\klconste$ we have $\xi\in \inter S_\AH$ for all such $\AH$.

On the other hand we consider the set $\mathcal{H}$ of all measures
on $\Zd$ with steplength $\step$, whose covariance matrices have smallest
eigenvalue greater or equal $\gamma$. $\mathcal{H}$ is obviously compact,
and the measures in $\mathcal{H}$ are all maximal (the covariance matrix of
a non maximal measure has zero as an eigenvalue). The function that maps
$(H,\xi)$ to the smallest eigenvalue of the covariance matrix of
$\AH_{t_\xi}$ is continuous in both arguments. Hence we find a constant
$\klconst\le\klconste$ such that the smallest eigenvalue of the covariance
matrix of $\AH_{t_\xi}$ is greater or equal $\gamma/2$ for all
$H\in\mathcal{H}$ and $\abs{\xi}\le\klconst$. This shows the claimed
independence.

Now we come back to the proof of the LCLT. We know
$\abs{\xi}=\abs{x/n}\le n^{-(4+\alpha)/9}$ for a fixed constant
$\alpha\in(0,1/2)$. Thus for
almost all
$n$ we have $\abs{\xi}\le\klconst$. It is sufficient to prove the
estimate for these $n$, since there is only a finite number of $n$ and
$x$ with
$\klconst n< \abs{x}\le n^{(5-\alpha)/9}$. We can cover these cases by
choosing $\const$ large enough. Moreover, using the compactness of
$\mathcal{H}$ again, we see that we find $\const$ depending on $d$,
$\step$, $\klew$ and $\alpha$ only.

So let $\xi \in S_\AG $ with $\abs{\xi}\le \klconst
\wedge n^{-(4+\alpha)/9}$. To estimate the first factor in \eqref{AG-n},
we use Taylor expansion for $I$ at zero.

Remember $I(\xi)= t_\xi\cdot \xi - \log Z(t_\xi)$. We
have $D I(\xi)=t_\xi$, and a simple computation yields $D^2 I(\xi)=
\V_\xi^{-1}$.
In general $D^k I(\zeta) (\xi,\dots,\xi)$ is a homogeneous
polynomial of degree $k$, whose coefficients depend (polynomially) on
$\V_\zeta^{-1}$ and on the first to
$k$th moments of $\AG_\zeta$.

Now let $T^{(k)}$ denote a homogeneous polynomial of degree $k$.
Using $I(0)=0$ and $D I(0)=0$, we obtain
\begin{align*}
I(\xi) &=  \xi \cdot (\V^{-1} \cdot \xi)/2 + T^{(3)}(\xi) +
T^{(4)}(\xi) + R_5(I,0,\xi).
\end{align*}

To deal with the
error term $R_5$ observe that, since the steplength of the tilted
measures is also bounded by $\step$, we can find upper bounds (depending
only on $d$ and $\step$) for their moments. Furthermore
we have $\abs{\V^{-1}\cdot\xi}\le\abs{\xi}/\klew$. Thus we know
\begin{align*}
\abs{R_5(I,0,\xi)}
&\le \text{const}(d,\step,\klew) \abs{\xi}^5,
\end{align*}
and hence
\begin{align}                                            \label{exp-nI}
\exp[-n I(\xi)] &= \exp[-n\,\xi \cdot (\V^{-1} \cdot \xi)/2]
\times \exp[-n T^{(3)}(\xi) -n T^{(4)}(\xi) -n
O((\xi)^{5}))]\notag\\[2mm]
   &= \exp[-n\,\xi \cdot (\V^{-1} \cdot \xi)/2]
\times \!\bigl[1-n T^{(3)}(\xi) -n T^{(4)}(\xi)+ n^2
T^{(6)}(\xi)\notag\\
&\qquad\qquad\qquad\qquad+\!\underbrace{O(n\abs{\xi}^{5}) +
\!O(n^2\abs{\xi}^{7}) +\! O( n^3\abs{\xi}^{9}) \exp[O( n \abs{\xi}^3)]}_{=
O(n^{-1-\alpha})
\text{, since } \abs{\xi}\le n^{-(4+\alpha)/9}}
\bigr]\notag\\ &= \exp[-x \cdot (\V^{-1} \cdot x)/(2n)]\\
&\quad\times \!\bigl[1-n^{-1/2} T^{(3)}(x/\sqrt{n}) -n^{-1} (T^{(4)}+
T^{(6)})(x/\sqrt{n})+O(n^{-1-\alpha})\bigr].\notag
\end{align}
The coefficients of the polynomials $T^{(k)}$ are depending polynomially
on $\V^{-1}$ and the moments of $\AG $ up to order four, and
$\abs{O(n^{-1-\alpha})}\le \text{const}(d,\step,\klew,\alpha)
\,n^{-1-\alpha}$.

As a next step we approximate the second factor in \eqref{AG-n},
$\AG _{t_\xi}^{*n}(x)$. This is done in the second part of the paper,
using Cram\'er-Edgeworth expansions for lattice
distributions. Applying the corollary on page \pageref{cor-bhatt} on the
measures $\AG _{t_\xi}$, we obtain
\begin{align*}
\babs{\AG _{t_\xi}^{*n}(x)
-&(2\pi n)^{-d/2} (\det \V_\xi)^{-1/2}\bigl[1
+ n^{-1} \con(\xi)\bigr]}
\le \cons \; n^{-(d+3)/2},
\end{align*}
where $\con(\xi)$ depends polynomially on the moments
of $\AG_{t_\xi}$ up to order four and on
the matrix $\V_\xi^{-1}$. Since $\klew/2$ is a lower bound for the
smallest eigenvalue of $\V_\xi$, the constant
$\cons$ can be chosen depending only on $d$, $\step$ and $\klew$.

We know that $\V_\xi$, $\V_\xi^{-1}$ and the other moments of
$\AG_{t_\xi}$ are analytic functions of $\xi$. Taylor expansion around
$\xi=0$ thus yields
\begin{align}                                       \label{Gxin}
\AG _{t_\xi}^{*n}(x)
&=(2\pi n)^{-d/2} (\det \V)^{-1/2}\notag\\
&\qquad\times\bigl[1
+ T^{(1)}(\xi) + T^{(2)}(\xi) + O(\xi^3) + n^{-1}[T^{(0)}+O(\xi)]
+O(n^{-3/2})\bigr]\notag\\
&=(2\pi n)^{-d/2} (\det \V)^{-1/2}\\
&\qquad\times\bigl[1
+ n^{-1/2} T^{(1)}(x/\sqrt{n}) + n^{-1} (T^{(0)}+T^{(2)})(x/\sqrt{n})
+O(n^{-1-\alpha})\bigr],\notag
\end{align}
with the same dependencies as in \eqref{exp-nI}.
Inserting \eqref{exp-nI} and \eqref{Gxin} in \eqref{AG-n} yields
\begin{align}                                  \notag 
\AG ^{*n}(x) &=
\varphi_{n\V}(x)\times
\!\bigl[1-n^{-1/2} P^{(3)}(x/\sqrt{n})
+n^{-1}P^{(6)}(x/\sqrt{n})+O(n^{-1-\alpha})\bigr],
\end{align}
where $P^{(k)}$ denotes a polynomial of order $k$. Here $P^{(3)}$
contains only odd order terms and $P^{(6)}$ only even order terms.

We still have to estimate $\varphi_{n\V}$. Since the coefficients of the
covariance matrix $\V$ are bounded by $\step^2$, we can bound the
maximal eigenvalue of $\V$ by $d\step^2$. Therefore
\begin{align*}
\varphi_{n\V} (x) &= (2\pi n)^{-d/2}(\det \V)^{-1/2}
\,\exp[- (x\cdot\V^{-1}x)/(2n)]\\
&\le (2\pi n \klew)^{-d/2} \,\exp[- x^2/(2nd\step^2)]\\
&\le (d\step^2/\klew)^{d/2}\,\varphi_{nd\step^2}(x)\\
&\le (2d\step^2/\klew)^{d/2}\,\varphi_{n 2d\step^2}(x).
\end{align*}
Choosing $\const$ large enough yields the claim for $\abs{x}\le
n^{(5-\alpha)/9}$.
\vspace{2mm}

{\em Case $\abs{x}> n^{(5-\alpha)/9}$:}

For \ane{big} $x$ we estimate $\AG ^{*n}(x)$ and
$\bigl[1+n^{-1/2}P_3 (x/\sqrt{n})
+ n^{-1}P_6 (x/\sqrt{n})\bigr]\varphi_{n\V}(x)$ separately.
For fixed natural $k$ we have $\abs{x}^k \exp [-x^2]\le \text{const}(k)
\exp[-x^2/\sqrt{2}]$ for all $x\in\Zd$.
Therefore  $\babs{1+n^{-1/2}P_3 (x/\sqrt{n})
+ n^{-1}P_6 (x/\sqrt{n})}\varphi_{n\V}(x)$ is bounded
by $\text{const}(d,\step,\klew)\, \varphi_{\sqrt{2}nd\step^2}(x)$. Now we
use the fact that there exists a natural number $k$ such that
$n^{1+\alpha}\le (x/\sqrt{n})^{k}$ ($k$ depending on $\alpha$). Thus
we find  $\const=\text{const}(d,\step,\klew,\alpha)$ such that for all $n$
\begin{equation}                                              \label{pol}
\babs{1+n^{-1/2}P_3 (x/\sqrt{n})
+ n^{-1}P_6 (x/\sqrt{n})}\varphi_{n\V}(x)\le
\const\, n^{-1-\alpha}
\varphi_{2nd\step^2}(x).
\end{equation}

Now consider $\AG ^{*n}(x)$. If $\xi\notin S_\AG $, $\AG ^{*n}(x)$
equals zero.
If $\xi\in \inter S_\AG $, we use Taylor expansion of $I(\xi)$ to obtain
\begin{align}                                    \label{Inachunten}
I(\xi) &= \int_0^1 (1-s) D^2I(s\xi)(\xi,\xi) ds
= \int_0^1 (1-s) (\xi\cdot\V_{s\xi}^{-1}\cdot\xi) ds \ge \xi^2/(2d\step^2).
\end{align}
In the last inequality we used that the maximal eigenvalue of $\V_{s\xi}$
is majorized by $d\step^2$ for all $\xi$ and $s$.
Inserting \eqref{Inachunten} into \eqref{AG-n} yields
$\AG ^{*n} (x) \le \exp[-x^2/ (2d \step^2 n)]$. Now the same arguments
that lead to \eqref{pol} yield the desired estimate for $\xi\in \inter
S_\AG$.

Finally we consider $\xi\in\partial S_\AG$. Approximating $\xi$ with a
sequence
$(\xi_k)_{k\in\N}$ of points in $\inter S_\AG$, the corresponding
sequence of measures $G_{t_{\xi_k}}$ converges to a measure $\AG'_\xi$
with support on a subset of the support of $G$, and $I(\xi_k)\to I(\xi)$.
Using \eqref{AG-n} and \eqref{Inachunten} we therefore obtain
\begin{align*}                                           
\AG ^{*n}(x) & = \exp[-n I(\xi)] \cdot \AG'_\xi{}^{*n}(x)
\le \exp[-x^2/(2d\step^2 n)].
\end{align*}
The rest of the argument proceeds as before.
\end{proof}

\section*{Cram\'er-Edgeworth expansions for lattice distributions}

This second part yields the bounds we use for the tilted measures in the
LCLT. The proof follows the arguments in Bhattacharya/Rao, \cite{BR}. The
statement here differs from their result in two points: On one hand, we
only want to approximate the walk at time $n$ up to order $n^{-(d+3)/2}$,
and we suppose bounded range of the measure. Thus we only look at a
special case of their result. On the other hand, we need more precise
control of the error terms. This means that most of the required lemmas
from \cite{BR} would have to be adjusted. Rather than doing that, we prove
everything we need directly.

The object of interest here is $\AH$, the single step distribution of
an aperiodic and maximal random walk on $\Zd$ with bounded steplength
$\step\in\N$.
   Let $\E$ denote the expectation and $\V$
the covariance matrix of $\AH$. The aim is
to give a good local bound for $\AH^{*n}(n\E)$, which we use in the proof
of the LCLT.

First we introduce some notations and basic
estimates (for more details see \cite{BR}, chapter 6). We denote by
$\huth$ the Fourier transform of $\AH$, that is
\begin{equation*}
\huth(t) \deff \sumx e^{it\cdot x} \,\AH(x).
\end{equation*}
Since $\AH$ has bounded steplength, all moments of $\AH$ exist.
We denote by $\mu_\nu$ the $\nu$th moment of $\AH$, where
$\nu$ is a $d$-dimensional integral vector, that is
\begin{equation*}
\mu_\nu \deff \sumx x^\nu \AH(x),\quad\text{and we have}\quad
i^{\abs{\nu}}\mu_\nu = (D^\nu \huth)(0).
\end{equation*}
The so called cumulants $\chi_\nu$ are given by
\begin{equation*}
i^{\abs{\nu}}\chi_\nu = (D^\nu \log\huth)(0)
\end{equation*}
and can be expressed in terms of moments by equating coefficients of
$t^\nu$ on both sides of the formal identity
\begin{equation}                                       \label{chisumme}
\sum_{\abs{\nu}\ge1} \chi_\nu \frac{(it)^\nu}{\nu!}
=\log\huth(t)
= \sum_{s=1}^\infty (-1)^{s+1} \; \frac{1}{s}\Bigl(\sum_{\abs{\nu}\ge1}
\mu_\nu \frac{(it)^\nu}{\nu!}\Bigr)^s.
\end{equation}

For $r\in\N$ we define the polynomial $\chi_r(z)$ in $z\in\C^d$ by
\begin{equation*}
\chi_r(z) = \sum_{\abs{\nu}=r} \frac{r!}{\nu!}\,\chi_\nu\, z^\nu.
\end{equation*}
We have $\chi_1(z)=\E\cdot z$ and $\chi_2(z)=z\cdot (\V \cdot z)$. For
$r\in\N$ we can bound
\begin{align}                                        \label{chiest}
\babbs{\frac{\chi_r(z)}{r!}}
&=\babbs{\sum_{\abs{\nu}=r} \chi_\nu
\frac{z^\nu}{\nu!}\,{}}
= \babbs{{}\,\sum_{k=1}^r (-1)^{k+1}
\frac{1}{k} \;
\sum_{\substack{\ssst\abs{\nu^{(1)}},\dots,\abs{\nu^{(k)}}\ge1\\
\sum_{j=1}^k \ssst\abs{\nu^{(j)}}=r}} \mu_{\nu^{(1)}}\dots\mu_{\nu^{(k)}}
\frac{z^{\sum_{j=1}^k\nu^{(j)}}}{\prod_{j=1}^k\nu^{(j)}!}\;{}}\notag\\
&\le \abs{z}^r \,\step^r\sum_{k=1}^r \frac{1}{k}
\Bigl(\underbrace{\sum_{\abs{\nu}\ge1}
\frac{1}{\nu!}}_{=e^d -1}\Bigr)^k
\le \abs{z}^r \: r\, [\step (e^d-1)]^r\notag\\
&\le r (\conste\abs{z})^r,
\end{align}
where $\conste\deff \step(e^d-1)$. Thus for
$\abs{t}<1/\conste$ the Taylor series $\sum_{r=1}^\infty
\frac{\chi_r(it)}{r!}$ is convergent, that is,
\begin{align}                                         \label{TElog}
\log \huth(t) &=\E\cdot it - \frac{t\cdot (\V \cdot t)}{2}
+\sum_{r=3}^\infty\frac{\chi_r(it)}{r!}.
\end{align}

\begin{lem}                                         \label{lem-bhatt}
There exists a constant $\cons$ such that for all $x$ in $\Zd$,
$\alpha\in (0,1/2)$ and for all  natural $n$,
\begin{align}                                       \label{diff}
\babbs{\AH^{*n}(x)
-&n^{-d/2} \Bigl[1
- \frac{\chi_3(D)}{6\,\sqrt{n}}
+ \frac{\chi_4(D)}{24\,n}
+ \frac{\chi_3^{*2}(D)}{72\,n}\Bigr]
\;\varphi_\V\Bigl(\frac{x-n\E}{\sqrt{n}}\Bigr)}\notag\\[2mm]
&\le \cons \; n^{-(d+3)/2},
\end{align}
where the formal polynomials $\chi_r(D)$ are given by
\begin{equation*}
\chi_r(D)f(x) \deff \sum_{\abs{\nu}=r} \frac{r!}{\nu!}\,\chi_\nu\, D^\nu
f(x)
\quad\text{and}\quad
\chi_r^{*2}(D)f(x) \deff (\chi_r(D)f*\chi_r(D)f)(x).
\end{equation*}
The constant $\cons$ depends only on the dimension $d$, the
steplength $\step$ of $\AH$ and on a lower bound for the smallest
eigenvalue of $\V$.
\end{lem}

\begin{proof}
First we use the Fourier inversion formula for aperiodic lattice measures
(see for example \cite{BR} (21.28)) to express $\AH$. We have
\begin{align}                                       \label{IF1}
\AH^{*n}(x)
&= (2\pi)^{-d}\int_{[-\pi,\pi]^d} \exp[-it\cdot x]\; \huth^n(t) \;dt
\notag\\
&= (2\pi\sqrt{n})^{-d}\int_{[-\sqrt{n}\pi,\sqrt{n}\pi]^d}
\exp[-it\cdot \frac{x}{\sqrt{n}} + n
\log\huth(\frac{t}{\sqrt{n}})]\;dt.
\end{align}
With the usual inversion formula for continuous measures we obtain
\begin{align}                                       \label{IF2}
&n^{-d/2} \Bigl[1
- \frac{\chi_3(D)}{6\,\sqrt{n}}
+ \frac{\chi_4(D)}{24\,n}
+ \frac{\chi_3^{*2}(D)}{72\,n}\Bigr]
\;\varphi_\V\Bigl(\frac{x-n\E}{\sqrt{n}}\Bigr) \\
= (2\pi&\sqrt{n})^{-d}
\int_{\Rd}\exp\bigl[-it\cdot\frac{x-n\E}{\sqrt{n}}-\frac{t\cdot(\V\cdot
t)}{2}\bigr]
\Bigl[1
+ \frac{\chi_3(it)}{6\,\sqrt{n}}
+ \frac{\chi_4(it)}{24\,n}
+ \frac{\chi_3^2(it)}{72\,n}\Bigr]
\,dt,\notag
\end{align}
where we used
\begin{equation*}
\widehat{D^\nu\varphi_\V}(t)=(-it)^\nu\widehat{\varphi_\V}(t)
=(-it)^\nu \exp\bigl[-\frac{t\cdot(\V\cdot t)}{2}\bigr].
\end{equation*}

Now let $\Nullum$ be a small ball around $0$ with radius $\varepsilon>0$
(precise conditions on $\varepsilon=\varepsilon(d,\step,\klew)$ will
be determined later). Using
\eqref{IF1} and
\eqref{IF2} we can bound the left hand side of
\eqref{diff} by the sum of three integrals $I$, $II$ and $III$ with
\begin{align}
I &= (2\pi\sqrt{n})^{-d}
\int_{\sqrt{n}\Nullum}\Babs{
\exp\bigl[n\log\huth(\frac{t}{\sqrt{n}})\bigr]
-\exp\bigl[it\cdot\frac{n\E}{\sqrt{n}}-\frac{t\cdot(\V\cdot t)}{2}\bigr]
\notag\\
&\hspace{6cm} \times
\Bigl[1 + \frac{\chi_3(it)}{6\,\sqrt{n}}
+ \frac{\chi_4(it)}{24\,n}
+ \frac{\chi_3^2(it)}{72\,n}\Bigr]}\,dt,
\notag\\
\intertext{and}
II &=
(2\pi\sqrt{n})^{-d}\int_{[-\sqrt{n}\pi,\sqrt{n}\pi]^d\setminus\sqrt{n}\Nullum}
\;\Babs{\huth^n(\frac{t}{\sqrt{n}})}\;dt,
\notag\\
\intertext{and}
III &= (2\pi\sqrt{n})^{-d}
\int_{\Rd\setminus \sqrt{n}\Nullum}
\exp\bigl[-\frac{t\cdot(\V\cdot t)}{2}\bigr]
\Babs{1 + \frac{\chi_3(it)}{6\,\sqrt{n}}
+ \frac{\chi_4(it)}{24\,n}
+ \frac{\chi_3^2(it)}{72\,n}}
\,dt.
\notag
\end{align}

We start with the essential term, the integral $I$. As a first condition
we let $\varepsilon\le 1/(2\conste)$. Then we have
$\conste \abs{t}/\sqrt{n}\le 1/2$ for all $t\in\sqrt{n}\Nullum$.  Denote
by $\klew$ (a lower bound for) the smallest eigenvalue of $\V$. We use
\eqref{TElog} and \eqref{chiest} to obtain
\begin{align*}
I &\le (2\pi\sqrt{n})^{-d}
\int_{\sqrt{n}\Nullum} \exp\bigl[-\frac{t\cdot(\V\cdot t)}{2}\bigr]
\notag\\
&\qquad\qquad\qquad\qquad \times
\Babs{
\exp\bigl[\sum_{r=3}^\infty n^{-\frac{r-2}{2}}\frac{\chi_r(it)}{r!}\bigr]
 - \Bigl[1 +
\frac{\chi_3(it)}{6\,\sqrt{n}} + \frac{\chi_4(it)}{24\,n}
+ \frac{\chi_3^2(it)}{72\,n}\Bigr]}\,dt
\notag\\
&\le (2\pi\sqrt{n})^{-d}
\int_{\sqrt{n}\Nullum} \exp\bigl[-\frac{\klew}{2}t^2\bigr]
\times \Bigl[\sum_{r=5}^\infty n^{-\frac{r-2}{2}}
\frac{\abs{\chi_r(it)}}{r!}
\notag\\
&\hspace{2.1cm}+ \frac{1}{2} \sum_{\substack{r,s\ge 3\\r+s\ge 7}}
n^{-\frac{r+s-4}{2}}
\frac{\abs{\chi_r(it)}\abs{\chi_s(it)}}{r!s!}
+\sum_{s=3}^\infty
\Bigl(\sum_{r=3}^\infty n^{-\frac{r-2}{2}}
\frac{\abs{\chi_r(it)}}{r!}
\Bigr)^s/s!\Bigr] \,dt
\notag\\
&\le (2\pi\sqrt{n})^{-d}
\int_{\sqrt{n}\Nullum} \exp\bigl[-\frac{\klew}{2}t^2\bigr]
\; \Bigl[n^{-3/2}\sum_{r=5}^\infty
r(\conste t)^r
+ \frac{1}{2} n^{-3/2}
\sum_{\substack{r,s\ge 3\\r+s\ge 7}}
rs(\conste t)^{r+s}\notag\\
&\hspace{8.4cm}
+\sum_{s=3}^\infty n^{-s/2}(2\conste t)^{3s}/s!\Bigr] \,dt
\notag\\
&\le n^{-(d+3)/2}\;(2\pi)^{-d}\int_{\Rd}
\exp[-\frac{\klew}{2} t^2]\,\text{const}(d,\step)\, \bigl[t^5 + t^7 +
t^9 \exp[\varepsilon (2\conste)^3) t^2]\bigr]\,dt
\notag\\
&\le n^{-(d+3)/2}\;(2\pi)^{-d}
\int_{\Rd}
\, \text{const}(d,\step,\klew) \;\exp[-\frac{\klew}{4} t^2]\,dt,
\notag
\end{align*}
if we choose $\varepsilon=\varepsilon(d,\step,\klew)$ small
enough. Therefore we have
\begin{align*}
I &\le \text{const}(d,\step,\klew) \; n^{-(d+3)/2}.
\end{align*}

Now we come to the error term II. Since $\abs{\huth}=1$ in zero and
$\abs{\huth}<1$ away from zero, there exists a constant $\constz<1$ with
$\abs{\huth(t)}\le\constz$ on $[-\pi,\pi]^d\setminus\Nullum$. As in the
proof of the LCLT we use the fact, that the set of all measures on $\Zd$
with steplength
$\step$ such that the smallest eigenvalue of the covariance matrix is
greater or equal $\klew$, is compact. Since the mapping $(H,t)\mapsto
\huth(t)$ is continuous in both arguments, the constant $\constz$ can be
chosen depending on $d$, $\step$ and $\klew$ only. Therefore we
have
\begin{align*}                                       
II &=
(2\pi)^{-d}\int_{[-\pi,\pi]^d\setminus\Nullum}
\;\babs{\huth^n(t)}\;dt
\le \constz^n \le \text{const}(d,\step,\klew)\; n^{-(d+3)/2}.
\end{align*}

$III$ is not much more complicated. Using \eqref{chiest} to majorize the
$\chi_r$, we can bound
\begin{align*}                                       
III &= (2\pi\sqrt{n})^{-d}
\int_{\Rd\setminus \sqrt{n}\Nullum}
\exp\bigl[-\frac{t\cdot(\V\cdot t)}{2}\bigr]
\Babs{1 + \frac{\chi_3(it)}{6\,\sqrt{n}}
+ \frac{\chi_4(it)}{24\,n}
+ \frac{\chi_3^2(it)}{72\,n}}
\,dt
\notag\\
&\le  (2\pi\sqrt{n})^{-d}
\int_{\Rd\setminus \sqrt{n}\Nullum}
\text{const}(d,\step)\exp\bigl[-\frac{\klew t^2}{3}\bigr]
\,dt
\notag\\
&\le \text{const}(d,\step,\klew)\, n^{-d/2}\exp\bigl[-\frac{\klew
\varepsilon^2 n}{4}\bigr]
\le \text{const}(d,\step,\klew)\; n^{-(d+3)/2}.
\end{align*}
The combination of these three estimates completes the proof of the Lemma.
\end{proof}

For the proof of the Local CLT we only need an approximation of
$\AH^{*n}(x)$ at the point $x=n\E$. In this case the approximating term
\begin{equation*}
n^{-d/2} \Bigl[1
- \frac{\chi_3(D)}{6\,\sqrt{n}}
+ \frac{\chi_4(D)}{24\,n}
+ \frac{\chi_3^{*2}(D)}{72\,n}\Bigr]
\;\varphi_\V\Bigl(\frac{x-n\E}{\sqrt{n}}\Bigr)
\end{equation*}
simplifies to
\begin{equation*}
n^{-d/2} \bigl[1
+ n^{-1}\bigl(\chi_4(D)/24 + \chi_3^{*2}(D)/72\bigr)\bigr]
\;\varphi_\V(0),
\end{equation*}
since the odd derivatives of a centered normal density vanish at zero.
Therefore we obtain the following corollary as a special case of the
lemma.
\begin{cor}                                         \label{cor-bhatt}
There exist constants $\cons$ and $\con$ such that for all natural $n$,
\begin{equation*}
\babs{\AH^{*n}(n\E)
-(2\pi n)^{-d/2} (\det \V)^{-1/2}\Bigl[1 + n^{-1} \con\bigr]}
\le \cons \; n^{-(d+3)/2}.
\end{equation*}
Here $\con$ is a constant depending polynomially on the moments
of $\AH$ up to order four (coming from the $\chi$ terms) and on the
matrix $\V^{-1}$ (coming from the derivatives of $\varphi_\V$).
The constant $\cons$ depends only on the dimension $d$, the
steplength $\step$ of $\AH$ and on a lower bound for the smallest
eigenvalue of $\V$.
\end{cor}



\end{document}